\documentstyle{amsppt}

\magnification=1200 \NoBlackBoxes \hsize=11.5cm \vsize=18.5cm

\def\M{\Cal M}

\def\Y{\Cal Y}

\def\C{\Bbb C}

\def\P{\Bbb P}

\def\L{\Cal L}
\def\O{\Cal O}

\def\y{\bar{y}}

\def\Sym{\text{Sym}}

\def\1/2{\frac{1}{2}}

\def\Y{\Cal Y}
\def\SS{\Cal S}

\def\2{{[2]}}
\def\l{\ell}
\topmatter
\title The (-1) twist \endtitle
\author
Ziv Ran\footnote{\raggedright{ Partially supported by NSA Grant
MDA904-02-1-0094} }
\endauthor

\address University of California, Riverside\endaddress
\email ziv\@math.ucr.edu\endemail \rightheadtext {-1}
\leftheadtext{} \abstract We prove some
 lower bounds on  certain twists
 of the canonical bundle of a
subvariety of a generic hypersurface in projective space. In
particular we prove that the generic sextic threefold contains no
rational or elliptic curves and no nondegenerate curves of genus
2.\endabstract

\endtopmatter\document
The purpose of this note is to extend the main result of [CR]
to cover the case of the $(-1)$ twist. The result is as follows.
\proclaim{Theorem }
Let $$X\in\L_d$$ be generic with
$$d(d-2)/8\geq 3n-1-k, d>n$$
 and
$$f:Y\to X$$ a
desingularization of an irreducible subvariety
of dimension $k$.
Set
$$t=\max(-1,-d+n+1+[\frac{n-k}{2}]).$$
Then  either
$$h^0(\omega_Y(t))> 0$$
or $f(Y)$ is contained in the union of the lines
lying on $X$.\par
In the case where $k=n-3$ we have furthermore that
if $h^0(\omega_Y(t))=0$,
then $f(Y)$ is ruled by lines.\endproclaim

The proof is largely identical to that in [CR], so we will
just indicate the differences. Assuming
$$h^0(\omega_Y(-1))=0,$$
we wish to show that $Y$ is contained in the
union of lines in $X$. Compared to [CR], display (6.3),
we may now assume only
$$-d+n+1+s\geq 0,$$
then Remark 6.3.1 allows us to assume that in general
$\l(y)$ is a 'bicontact line' to $X_{F(y)}$, i.e. meets it
in 2 points with multiplicities $r, d-r$ for some
$1\leq r\leq d-1$. This gives rise to a map
$$\Y\to\Delta_{r,d-r}$$
 where
$\Delta_{r,s}$ is the closure of
$$\{(\l,x,x',F):x\neq x', \l.X_F\geq rx+sx'\}
\subset G\times\P^n\times\P^n\times S.$$
Note that $\Delta_{r,s}$ is an open subset of
a vector bundle over the double incidence variety
$$I^2:=\{(\l,x,x'):x,x'\in\l\},$$
and in particular is smooth.
Below we shall compute the canonical bundle of
$\Delta_{r,s}$. Once this is done, the proof may be concluded
by copying the arguments of [CR], beginning with Lemma 6.4.
\par
We shall assume
henceforth that
$$r\geq s.$$
Consider the product
$$\Delta_r\times\P^n,$$
with line bundle
$$\O_2(-1):=p_2^*(\O(-1)).$$
Then the zero-scheme $\Delta_{r,0}$ of the natural map
$$\O_2(-1)\to Q$$
is just the locus
$$\{(\l,x,x',F):(\l,x,F)\in \Delta_r, x'\in\l\}$$
and admits a natural map to $I^2.$

Clearly, we have
$$\omega_{\Delta_{r,0}}=\omega_{\Delta_r}\otimes\O_G(1)
\otimes\O_2(-n+1).$$ Now consider on on $\Delta_{r,0}$ the
zero-scheme of the natural map
$$\O_S(-1)\to \O_2(d).$$
This consists of $\Delta_{r,1}$ plus the pullback of the diagonal
divisor
$$D=\{(\l,x,x)\}\subset I^2.$$
Since $D$ is the degeneracy locus of the natural map
$$\O_1(-1)\oplus\O_2(-1)\to \SS,$$
it is easy to see that
$$\O(D)=\O_G(-1)\otimes\O_1(1)\otimes\O_2(1),$$
hence by the adjunction formula we have
$$\omega_{\Delta_{r,1}}=\omega_{\Delta_{r,0}}\otimes
\O_G(1)\otimes\O_1(-1)\otimes\O_2(d-1).$$ Now we can argue as
above, considering the natural injection
$$\O_2(-1)\to \SS$$
and the induced filtration $F^._2$ on $\Sym^d(\SS^{\text{v}}).$
The zero scheme of the natural map
$$\O_S(-1)\to F^1_2/F^{s}_2$$
consists of $\Delta_{r,s}$ plus $D$, i.e. $\Delta_{r,s}$ is a zero
scheme of $(F^1_2/F^s_2)(-D)$, and therefore
$$\omega_{\Delta_{r,s}}=\omega_{\Delta_{r,1}}
\otimes\det(F^1_2/F^s_2)\otimes\O(-(s-1)D).$$ This finally gives
the formula (valid for $r\geq s$)
$$\omega_{\Delta_{r,s}}=$$
$$\O_G(\frac{r(r-1)+s(s-1)}{2}
+s+1-n)\otimes\O_1(r(d-r)-n+r-s+1)\otimes\O_2 (s(d-s)-n+1).$$
\remark{Remark} It is clear from [CR, Remark 6.3.1] that to treat
the case of the $(-1)$ twist would require dealing with
the bicontact locus. Shortly after Herb Clemens and
the author observed this, and while they were trying
to establish the $(-1)$ case, the author was informed by Clemens
that Pacienza, who had earlier requested and obtained from
Clemens a working draft of [CR], was claiming to be able to do
the $(-1)$ case by computing the canonical bundle of
$\Delta_{r,d-r}$, which is a fairly obvious extension of
the method of [CR]. The author feels strongly that the $(-1)$
case should have been incorporated into [CR], but was unable
to convince Herb Clemens to agree.\endremark

\heading Appendix\endheading
The purpose of this appendix is to give an alternate and
somewhat shorter proof of the main Lemma  (Lemma 1.2) of [CLR],
which goes as follows.
\proclaim{Lemma } Let $Y\subset\P (\L)$
be an irreducible subvariety spanning a $\P^{p+1}.$  Fix integers
$r,k$ with $0\le r-1\le k\le d$. Let $\Cal L_d\to \C^N$ be a
linear map of vector spaces such that for $y_1,\dots,y_{r-1}$
general points of $Y$, the restriction to $\Cal
L_d(-y_1-\dots-y_{r-1})$ induces a surjection
$$\Psi:\Cal
L_d(-y_1-\cdots-y_{r-1})\to\C^{k+1}.$$
Then for a general choice
of elements $h_{k+1},\dots, h_d\in\Cal L$ and for general subsets
$Y_1,\dots,Y_k \subset Y$ each of cardinality $p$, with $Y_i\ni
y_i$, $i=1,\dots,r-1$,  the restriction of $\Psi$ to the subspace
$\Cal L(-Y_1)\cdots \Cal L(-Y_k)h_{k+1}\cdots h_d$ surjects.
\endproclaim
\demo{proof} By a {\it {good chain}} in $\L$ we mean a (connected)
chain whose components are straight lines (i.e. pencils) of the
form $\L(-S)$ where $S$ is a general $p$-tuple in $Y$ and whose
'vertices' (i.e. singular points) are general in $\P (\L )$. Clearly two
general elements of $\L$ can be joined by a good chain.
\def\M{\Cal M}
Let $\M_d\subset\L_d$ be the set of monomials. By a {\it {good
chain}} in $\M_d$ we mean a chain which is a union of subchains of
the form $\Cal C_i'=h_1\cdots\Cal C_i\cdots h_d$ where $\Cal C_i$
is a good chain in $\L$. It is easy to see that two general
monomials can be joined by a good chain. \par Next, let us say
that a monomial $h_{d-e+1}\cdots h_d\in\M_e$ is rel $\bar{y},
\bar{y}=\{y_1,...,y_{r-1}\}$ if $h_i\in\L(-y_i), i=d-e+1,...,r-1$
(this condition is vacuous if $r-1<d-e+1$); denote by
$\M_e(-\bar{y})$ the set of these. Again it is easy to see that
two general elements of $\M_e(-\bar{y})$ can be joined by a good
chain within  $\M_e(-\bar{y})$.\par Now to prove the Lemma it
suffices to show by induction on $q, 0\le q\leq k$, that, with the
above notations,

$$ \dim\Psi (\L(-Y_1)\cdots\L(-Y_q)h_{q+1}\cdots h_d)\geq\min
(q+1, k+1)$$

(for general choices rel $\bar{y}$). For $q=0$ this is clear as
$\Psi$ is nonzero on a general element of $\M_d(-\bar{y})$,
because these span $\L_d(-\bar{y}$). Assume it is true for $q$ and
false for $q+1$, and suppose first that $q+1\leq r-1$. Now because
$$Z:=\Psi (\L(-Y_1)\cdots\L(-Y_q)h_{q+1}\cdots h_d) =\Psi
(\L(-Y_1)\cdots\L(-Y_q)\L(-Y_{q+1})h_{q+2}\cdots h_d)$$ by
assumption (i.e $Z$ doesn't move as $h_{q+1}$ varies in
$\L(-Y_{q+1})$, and because two general elements of $\L(-y_{q+1})$
can be joined by a good chain (whose components are of the form
$\L(-Y_{q+1})$), it follows that $Z$ is independent of
$h_{q+1}\in\L(-y_{q+1})$, fixing the other $h$'s. Since
$y_1,...y_{r-1}$ are all interchangeable, it follows that a
similar statement holds for any permutation of them. In particular
$Z$ contains all pencils of the form $\Psi
(h_1\cdots\L(-Y_j)\cdots h_d ), j=1,...,r-1$ and from connectivity
of $\M_{r-1}(-\y)$ by good chains it follows that in fact $Z =\Psi
(\L_{r-1}(-\y)h_r\cdots h_d)$. As for the remaining $h_j$'s, say
$j=r$, we may pick $y_r\in f^*(h_r)_0$ and apply similar reasoning
to $y_1,...,y_{r-1},y_r$ in place of $y_1,...,y_{r-1}$ . The
foregoing argument yields that $Z$ is independent of
$h_r\in\L(-y_r)$, fixing $y_r$ and the other $h$'s (indeed that
$Z=\Psi (L_r(-\y-y_r)h_{r+1}\cdots h_d)$. Then another similar
argument  with good chains of $h_j$'s not fixing $y_j$ yields
easily that  $Z$ is actually independent of $h_j\in\L $. Now since
we can connect two general elements of $\M_d(-\bar{y})$ by a good
chain, we conclude that  $Z=\Psi (\L_d (-\bar {y}))$, which is a
contradiction.

The case $q+1\geq r$ is similar but simpler: we may conclude
directly that $Z$ is independent of $h_{q+1}\in\L(-Y_{q+1})$,
hence of $h_j\in\L(-Y_j) $ for all $r\leq j\leq k$ and use good
chains as above to deduce a contradiction.\qed
\enddemo

\enddocument